\documentclass[11pt]{amsart}
\usepackage{xcolor}
\usepackage{hyperref}
\hypersetup{
    colorlinks = false,
    linkbordercolor = {white},
    citebordercolor = {white}
    }
\usepackage[capitalise,noabbrev]{cleveref}
\usepackage[all]{xy}

\usepackage{amsmath,amscd}
\usepackage{amssymb}
\usepackage{mathtools}
\usepackage{stmaryrd}
\usepackage{url}
\usepackage{tikz-cd}
\usepackage{enumitem}
\usepackage{fullpage}
\usepackage{musicography}
\usepackage{physics}

\usepackage{tikz}

\usepackage{enumitem,kantlipsum}

\author{Prashanth Sridhar}

\newcommand{\Addresses}{{
	\vskip\baselineskip
  	\footnotesize
  	\noindent %\textsc{Department of Mathematics, Purdue University} 
    \textsc{Department of Mathematics, University of Alabama, Tuscaloosa, AL, USA} \par\nopagebreak
	\noindent \textit{E-mail addresses:} \texttt{psridhar1@ua.edu}
 }}

\numberwithin{equation}{section}

\newtheorem{lem}[equation]{Lemma}
\newtheorem{theorem}[equation]{Theorem}

\newtheorem{cor}[equation]{Corollary}

\newtheorem{claim*}{Claim}

\theoremstyle{definition}

\newtheorem{setup}[equation]{Setup}

\theoremstyle{remark}

\newtheorem{rem}[equation]{Remark}

\newcommand{\mfrak}[1]{\mathfrak{#1}}

\newcommand{\m}{\mfrak{m}}
\newcommand{\n}{\mfrak{n}}

\newcommand{\projdim}{\operatorname{proj\,dim}}

\newcommand{\depth}{\operatorname{depth}}

\newcommand{\Hom}{\operatorname{Hom}}

\newcommand{\Mod}{\operatorname{Mod}}

\def\nc{\newcommand}
\nc{\on}{\operatorname}
\nc{\bideg}{\on{bideg}}
\nc{\xra}{\xrightarrow}
\def\phi{\varphi}

\nc{\into}{\hookrightarrow}
\nc{\onto}{\twoheadrightarrow}
\nc{\LL}{\mathbf{L}}
\nc{\RR}{\mathbf{R}}
\nc{\Perf}{\on{Perf}}
\nc{\nat}{\natural}
\nc{\tors}{\on{tors}}
\nc{\Tors}{\on{Tors}}

\def\Mod{\on{Mod}}
\nc{\qgr}{\on{qgr}}
\nc{\Qgr}{\on{Qgr}}
\nc{\fQgr}{\on{Qgr}^{\on{f}}}
\nc{\colim}{\on{colim}}

\nc{\Ext}{\on{Ext}}

\nc{\om}{\omega}
\nc{\w}{\widetilde}
\nc{\PP}{\mathbb{P}}
\nc{\mf}{\on{mf}}
\nc{\OO}{\mathcal{O}}
\nc{\Proj}{\on{Proj}}
\nc{\Qcoh}{\on{Qcoh}}
\nc{\coh}{\on{coh}}
\nc{\Tor}{\on{Tor}}
\nc{\Modf}{\Mod^{\on{f}}}

\nc{\ce}{\coloneqq}

\nc{\Com}{\on{Com}}
\nc{\A}{\mathcal{A}}
\nc{\B}{\mathcal{B}}
\nc{\C}{\mathcal{C}}
\nc{\Sh}{\on{Sh}}
\nc{\QCoh}{\on{QCoh}}
\nc{\Coh}{\on{Coh}}
\nc{\fQCoh}{\QCoh^{\on{f}}}
\nc{\ov}{\overline}
\nc{\End}{\on{\underline{End}}}
\def\MR#1{}

\nc{\Qgrf}{\Qgr^{\on{f}}}
\nc{\uHom}{\underline{\Hom}}
\nc{\Inj}{\mathrm{Inj}}
\nc{\proj}{\mathrm{Proj}}
\nc{\spec}{\mathrm{Spec}}

\def\A{\mathcal{A}}

\usepackage{microtype}
\usepackage{comment}

\begin{document}
\title{A note concerning the vanishing of local cohomology for roots in mixed characteristic}

\begin{abstract}
The goal of this note is to record the following curious fact: let $(S,\n)$ be an unramified regular local ring of mixed characteristic $p>0$ and dimension $d$. Let $L$ denote the quotient field of $S$ and $K=L(\omega)$ with $\omega^p\in L$. Let $R$ denote the integral closure of $S$ in $K$. Then $R$ is Cohen-Macaulay if and only if $\mathrm{H}^{d-1}_{\n}(R)=0$, i.e., the obstruction to the Cohen-Macaulayness of $R$ lies in a single local cohomology module. Furthermore, this is equivalent to the dual module $\Hom_S(R,S)$ satisfying Serre's condition $(S_3)$.
\end{abstract}


\numberwithin{equation}{section}

\maketitle
\setcounter{tocdepth}{1}
%\tableofcontents

\setcounter{section}{1}

%\section{ABCD}

This paper proves the following curious fact:

\begin{theorem}\label{thm:main}
    Let $(S,\n)$ be an unramified regular local ring of mixed characteristic $p>0$ and dimension $d$. Let $L$ denote the quotient field of $S$ and $K=L(\omega)$ with $\omega^p\in L$. Let $R$ denote the integral closure of $S$ in $K$. Then the following are equivalent:

    \begin{enumerate}
        \item $R$ is Cohen-Macaulay.

        \item $\mathrm{H}^{d-1}_{\n}(R)=0$.

        \item $R^*:=\Hom_S(R,S)$ satisfies $(S_3)$ as an $S$-module.

        \item $R^*:=\Hom_S(R,S)$ satisfies the $(S_3)$ condition as an $S$-module for all $Q\in V(p)$.
    \end{enumerate}
\end{theorem}

As an immediate corollary of \Cref{thm:main} and \cite[Theorem 3.8]{Evans_Griffith}, which is attributed to \cite{stable_module_theory}, we have

\begin{cor}
    With notation as in \Cref{thm:main}, the following are equivalent:

    \begin{enumerate}
        \item $R$ is an $n$-th syzygy over $S$ for all $1\leq n\leq d$.

        \item $R^*:=\Hom_S(R,S)$ is a third syzygy over $S$.

        \item $R^*:=\Hom_S(R,S)$ is an $n$-th syzygy over $S$ for all $1\leq n\leq d$.
    \end{enumerate}
\end{cor}

\Cref{thm:main} concerns the simplest modular Kummer extension of an unramified regular local ring. This appears in several contexts: as an example of the failure of a result of Roberts in \cite{RO} in the modular case, see \cite{KO,DK_99,GRIFFITH2015502,katz_sridhar_roberts}; in the context of existence of small Cohen-Macaulay modules in \cite{DK_99,SRIDHAR2021100,PS,DK}; and in applications to singularity theory in \cite{benozzo2025boundspluspurethresholdshypersurfaces}. We point out in \Cref{rem:syzygy_theorem} why the phenomenon in \Cref{thm:main} can be viewed as a mirror image result in this context to the syzygy theorem of Evans and Griffith, \cite[Theorem 1.1]{syzygy_problem_evans_griffith}, which is now known to be true in mixed characteristic \cite{andre}.

We recall that a module $M$ over a commutative ring $A$ satisfies $(S_n)$ if for all $P\in \mathrm{Spec}(A)$, $\mathrm{depth}(M_P)\geq \mathrm{min}\{n,\dim(M_P)\}$. If $A$ is a domain with field of fractions $K$, and normalization $R$, the ideal $\{x\in K\:|\: xR\subseteq A\}\subseteq A$ is called the conductor ideal of $A$. It can be canonically identified with the $A$-module $\Hom_A(R,A)$, see for example \cite[Lemma 2.4.2]{Huneke_Swanson}. 

\begin{proof}[Proof of Theorem 1.1]
    The implication (3) $\implies$ (4) is clear, while (4) $\implies$ (3) follows by working over $S[1/p]$ and applying \cite[Proposition 2.4]{SRIDHAR2021100}. The implication (1) $\implies$ (2) is also clear. To finish the proof, we will show (2) $\implies (3)$ and $(3)$ $\implies$ (1).

For the remainder of the proof, we will assume notation as specified below:

\begin{setup}
  We can and will assume $\omega^p=f\in S$. Let $a_1,\dots,a_n\in S$ be the divisors of $f$ in $S$ that satisfy $a_i^2|f$ (there could be none). Let $A$ be the hypersurface $S[X]/(X^p-f)\simeq S[\omega]\subseteq K$, where $X$ is an indeterminate over $S$. Let $\mathfrak{C}$ denote the conductor ideal of $A$.
\end{setup}

For a fractional ideal $I\subseteq K$, let $I^{-1}$ denote the fractional ideal $\{x\in K\:|\:xI\subseteq A\}$. We now record some standard facts in the lemma below.

\begin{lem}\label{lem}
    Assuming the setup above, the following hold.

\begin{enumerate}
    \item $R=\mathfrak{C}^{-1}\simeq \Hom_A(\mathfrak{C},A)$ and $\mathfrak{C}=R^{-1}\simeq \Hom_A(R,A)$, where the isomorphisms are as $A$-modules.
    \item $\mathfrak{C}$ satisfies $(S_2)$ as an $A$-module and is an unmixed ideal of height one.
    \item $R$ is Cohen-Macaulay if and only if $A/\mathfrak{C}$ is Cohen-Macaulay.
    \item $\Hom_A(R,A)\simeq \Hom_S(R,S)$ as $S$-modules.
\end{enumerate}
\end{lem}
\begin{proof}[Proof of Lemma 1.4] (1) By definition, $\mathfrak{C}=R^{-1}$. Moreover, the map $R^{-1}\to \Hom_A(R,A)$, $x\mapsto (r\mapsto r\cdot x)$ is an isomorphism of $A$-modules. Since the $A$-modules $R$ and $\mathfrak{C}^{-1}$ are birational, satisfy Serre's condition $(S_2)$ and agree in codimension one ($A$ is Gorenstein), they agree as $A$-submodules of $K$, see for instance \cite[15.10(i)]{Huneke_Swanson}. Finally, $R$ is a reflexive $A$-module since $A$ is a $(S_2)$-ring, $R$ is an $(S_2)$ $A$-module and is reflexive in codimension one over $A$ (in codimension one, being torsion free is the same as being reflexive over a Gorenstein ring), see \cite[Proposition 1.4.1]{Bruns-Herzog}. Thus, the isomorphism $\mathfrak{C}\simeq \Hom_A(R,A)$ implies $\Hom_A(\mathfrak{C},A)\simeq R$.
\\
(2) Since $A$ is $(S_2)$ and (1) says $\mathfrak{C}$ is an $A$-dual, $\mathfrak{C}$ is $(S_2)$ over $A$ and is hence unmixed of height one, see for instance \cite[15.10(i)]{Huneke_Swanson}.
\\
(3) See \cite[Proposition 2.11]{SRIDHAR2021100}.
\\
(4) See for instance \cite[Remark 4.4]{SRIDHAR2021100}).
\end{proof}

    \par We begin by considering $(3)$ $\implies$ (1). By hypothesis and \Cref{lem}(4) and \Cref{lem}(1), $\mathfrak{C}$ satisfies $(S_3)$ as an $S$-module. Since $S\subseteq A$ is module finite, $\mathfrak{C}$ satisfies $(S_3)$ as an $A$-module as well. Considering the short exact sequence 

    \begin{equation}\label{eqn:natural_cond}
  0\to \mathfrak{C} \to A \to A/\mathfrak{C} \to 0
\end{equation}

one sees that $A/\mathfrak{C}$ satisfies $(S_2)$ as an $A$-module, equivalently $(S_2)$ as a ring.

We may assume $A$ is not normal or equivalently that $\mathfrak{C}$ is not the unit ideal since in this case the theorem holds vacuously. Set $F(X):=X^p-f$. Since $K/L$ is separable, $F'(\omega)=p\cdot\omega^{p-1}\in \mathfrak{C}$, see \cite[Theorem 12.1.1]{Huneke_Swanson}. This also implies $p\cdot f\in \mathfrak{C}$. Since $\mathfrak{C}$ is unmixed, the associated primes of $\mathfrak{C}$ are amongst those of $p$ or $f$. Now suppose any of the following hold:

\begin{enumerate}
    \item $f$ is square free, i.e., there does not exist a divisor $a$ of $f$ in $S$ such that $a^2|f$.

    \item $\mathfrak{C}$ does not have a primary component containing $p$.
\end{enumerate}

In the case of (1), by \cite[Lemma 3.2]{DK_99}, $R$ is Cohen-Macaulay. In the case of (2), by \cite[Lemma 3.3 and 3.4]{DK_99}, $R$ is Cohen-Macaulay. Thus, the theorem holds vacuously in either of these cases. Hence, we may and will assume that $f$ is not square free and that $\mathfrak{C}$ has a primary component containing $p$.

\par Now let $P$ denote the primary component of $\mathfrak{C}$ containing $p$ (there is a unique one, see the end of section $2$ in \cite{DK_99}) and $J_i$ the primary component of $\mathfrak{C}$ corresponding to $a_i$, so that $\mathfrak{C}=P\cap J_1\cap \dots \cap J_n$. Consider the natural short exact sequence of $A$-modules

\begin{equation}
\label{eqn:natural}
  0\to A/\mathfrak{C} \to A/P \bigoplus A/(J_1\cap \dots \cap J_n) \to A/(P+J_1\cap \dots \cap J_n) \to 0
\end{equation}

Now we note that the proof of \cite[Proposition 2.11]{SRIDHAR2021100}, shows that if $I\subseteq A$ is a height one unmixed ideal, then $I^{-1}$ is a maximal Cohen-Macaulay $A$-module if and only if $A/I$ is a Cohen-Macaulay ring. Combining this with the proof of \cite[Lemma 3.2]{DK_99} shows that $A/P$ is Cohen-Macaulay. Moreover, \cite[Lemma 3.3 and Lemma 3.4]{DK_99} imply that $A/(J_1\cap \dots \cap J_n)$ is Cohen-Macaulay. Thus, by \Cref{lem}(3), $R$ is Cohen-Macaulay if and only if $\depth(A/(P+J_1\cap \dots J_n))\geq \depth(A)-2$. Another application of \cite[Lemma 3.2, Lemma 3.3 and Lemma 3.4]{DK_99} shows that $P+J_1\cap \dots J_n$ is of the form $(p,\omega-h, \omega^{g_s}, \omega^{g_{s-1}}c_1,\dots, \omega^{g_1}c_{s-1},c_s)$ for suitable non-units $h,c_1,\dots,c_s\in S$ and integers $g_1,\dots, g_s$. Now note that $X^p-f\in (p,X-h)$, see for example the last paragraph of Section 2 in \cite{DK_99}. Setting 
\[
\widetilde{J}:=(p, X-h, X^{g_s}, X^{g_{s-1}}c_1,\dots, X^{g_1}c_{s-1},c_s)\subseteq S[X]\]

we have $A/(P+J_1\cap \dots J_n)\simeq S[X]/(\tilde{J})\simeq \overline{S}/\overline{I}$ for $\overline{S}=S/pS$ and
$\overline{I}=(h^{g_s}, h^{g_{s-1}}c_1, \dots, h^{g_1}c_{s-1},c_s)\overline{S}$. Thus, $R$ is Cohen-Macaulay if and only if $\depth(\overline{S}/\overline{I})\geq \depth(\overline{S})-1$ if and only if $\projdim_{\overline{S}}(\overline{S}/\overline{I})\leq 1$ if and only if $\overline{I}$ is a free $\overline{S}$-module.

Since $A/\mathfrak{C}$ is $(S_2)$, it follows from \Cref{eqn:natural} that for all $Q\in \mathrm{Spec}(A)$ with $\mathrm{ht}(Q)\geq 3$, $\mathrm{depth}_{A_Q}(A/(P+J_1\cap \dots \cap J_n))_Q\geq 1$; equivalently, for all $Q\in \mathrm{Spec}(\overline{S})$ with $\mathrm{ht}(Q)\geq 2$, $\mathrm{depth}_{\overline{S}_Q}(\overline{S}/\overline{I})_Q\geq 1$. It then follows from the canonical short exact sequence

 \begin{equation}
  0\to \overline{I} \to \overline{S} \to \overline{S}/\overline{I} \to 0
\end{equation}

that $\overline{I}$ satisfies $(S_2)$ as an $\overline{S}$-module. Since $\overline{I}$ is torsion free, $\overline{I}$ is reflexive in codimension one over $\overline{S}$. Finally, since $\overline{S}$ satisfies $(S_2)$, it follows from \cite[Proposition 1.4.1]{Bruns-Herzog} that $\bar{I}$ is a reflexive $\bar{S}$-module. Since $\bar{I}$ has rank one (it is a submodule of the field of fractions of $\bar{S}$) and $\bar{S}$ is factorial, $\bar{I}$ is free, see for instance \cite[Theorem 5.3]{Samuel}. Hence $R$ is Cohen-Macaulay and this completes the proof of (3) $\implies$ (1). 

\par Finally, we consider (2) $\implies$ (3). We need to show $\Hom_S(R,S)\simeq \mathfrak{C}$ (\Cref{lem}(4) and \Cref{lem}(1)) satisfies $(S_3)$ as an $S$-module. Since $S\subseteq A$ is module finite, using the short exact sequence \ref{eqn:natural_cond}, it suffices to show $A/\mathfrak{C}$ satisfies $(S_2)$ as an $A$-module. Now note that for any $0\neq x\in \mathfrak{C}$, since $R=\mathfrak{C}^{-1}$ (\Cref{lem}(1)), $R$ can be identified with the $A$-submodule $(xA:_A\mathfrak{C})/x$ of $K$. In particular, $R\simeq (xA:_A\mathfrak{C})$ as $A$-modules. Moreover, since $\mathfrak{C}$ is unmixed of height one (\Cref{lem}(2)), it follows that $\mathfrak{C}$ and $(xA:_A\mathfrak{C})$ are linked, see for example \cite[Proposition 2.2]{schenzel}. Noting that $A$ is local with maximal ideal $\m=(\n,\omega-h)$, one sees from the long exact sequence in local cohomology associated to 

\begin{equation}
  0\to (xA:_A\mathfrak{C}) \to A \to A/(xA:_A\mathfrak{C}) \to 0
\end{equation}

that $\mathrm{H}^{d-2}_{\m}(A/(xA:_A\mathfrak{C}))=0$. Applying \cite[Theorem 4.1]{schenzel}, we then see that $A/\mathfrak{C}$ satisfies $(S_2)$ as an $A$-module and the proof is complete.
\end{proof}

\begin{rem}\label{rem:syzygy_theorem}
The syzygy theorem of Evans and Griffith, \cite[Theorem 1.1]{syzygy_problem_evans_griffith} says that any module with full support over a regular local ring that satisfies Serre's condition $(S_k)$ and has rank strictly less than $k$ must be free. This was originally stated in the equicharacteristic case; since the direct summand theorem \cite{andre} implies this result, it is now known to be true in the mixed characteristic setting too. In particular, if $S$ is a regular local ring, $L$ its field of fractions and $K/L$ a finite field extension of degree $k$, then in view of the direct summand theorem (and the syzygy theorem), the integral closure of $S$ in $K$ is Cohen-Macaulay if it satisfies $(S_k)$. Arguing as in the proof of \Cref{thm:main} and applying \cite[Theorem 4.1]{schenzel}, the requirement that $R$ satisfies $(S_k)$ can be reinterpreted as requiring the dual module $\Hom_S(R,S)$ to have its $k-2$ sub-maximal local cohomology modules vanish. In the special context of this paper, \Cref{thm:main} provides a mirror sufficient condition for $S$-freeness of $R$, interchanging the roles of $R$ and its $S$-dual.
    
\end{rem}

\begin{rem}
    For an explicit example of $R$ as in \Cref{thm:main} such that $R$ is not Cohen-Macaulay, we refer the reader to \cite[Example 3.10]{DK_99}.
\end{rem}

\subsection*{Acknowledgments} I thank Linquan Ma for helpful conversations. I also thank the referee for their helpful comments and careful reading of the paper.

\bibliographystyle{amsalpha}
\bibliography{references}
\Addresses
\end{document}